\documentclass[a4paper,12pt]{article}
\usepackage{amsfonts,amsthm,amsmath,amssymb,graphicx}

\title{On some non-conformal fractals}

\date{}

\author{Micha\l \ Rams \thanks{supported by EU FP6
    Marie Curie programme SPADE 2 and by Polish MNiSW Grant NN201 0222 33 `Chaos, fraktale i dynamika konforemna'}\\
\normalsize Institute of Mathematics, Polish Academy of Sciences\\
\normalsize ul. \'Sniadeckich 8, 00-950 Warszawa, Poland\\
\normalsize e-mail: rams@impan.gov.pl}

\theoremstyle{plain}
\newtheorem{lem}{Lemma}[section]
\newtheorem{prop}[lem]{Proposition}
\newtheorem{thm}[lem]{Theorem}

\theoremstyle{definition}

\theoremstyle{remark}
\newtheorem*{rem}{Remark}

\numberwithin{equation}{section}

\newcommand{\N}{\mathbb N}

\renewcommand{\epsilon}{\varepsilon}

\begin{document}

\maketitle

\def\thefootnote{}
\footnote{1991 {\it Mathematics Subject Classification}: 28A78,
28A80}
\def\thefootnote{\arabic{footnote}}
\setcounter{footnote}{0}

\begin{abstract}
This paper presents a simple method of calculating the Hausdorff
dimension for a class of non-conformal fractals.
\end{abstract}

\section{Introduction}

An iterated function scheme acting on a complete metric space $X$
is a finite family of contracting maps ${\cal F}=\{f_k\}_{k=1}^n;
f_k:X\to X$. As noted by Hutchinson \cite{Hu}, the related
multimap

\[
F(\cdot)=\bigcup_{k=1}^n f_k(\cdot)
\]
(acting on the space $B(X)$ of nonempty compact subsets of $X$,
considered with the Hausdorff metric) is also a contraction.
Hutchinson proved that if $X$ is complete, so is $B(X)$. Hence, by
the Banach fixed point theorem, there exists a unique nonempty
compact set $\Lambda$ satisfying

\[
\Lambda=F(\Lambda)=\lim_{n\to\infty} F^n(A).
\]
The limit does not depend on the choice of $A\in B(X)$. $\Lambda$
is called the limit set of the iterated function scheme $\cal{F}$.

By similar reasoning, if we have a finite number of iterated
function schemes $\{{\cal F}_i\}_{i=1}^m$ acting on $X$ and apply
them in any order, the pointwise limit

\[
\Lambda_\omega = \lim_{n\to\infty} F_{\omega_1}\circ \ldots \circ
F_{\omega_n}(A)
\]
exists for all $\omega\in \Omega= \{1,\ldots,m\}^\N$ and does not
depend on $A\in B(X)$.

The question we want to answer (motivated by \cite{Lu}, see also
\cite{N}, \cite{GL}, \cite{GL2} and the incoming paper \cite{Re})
is: when the iterated function schemes ${\cal F}_i$ are of some
special class (for which we can calculate the Hausdorff dimension
of the limit set of any deterministic iterated function scheme
from this class) and the sequence $\omega$ is chosen, what will be
the value of the Hausdorff dimension of $\Lambda_\omega$?

We will present a simple method of dealing with this question,
working for Lalley-Gatzouras maps \cite{LG}, Bara\'nski maps
\cite{B} and higher dimensional affine-invariant sets of Kenyon
and Peres \cite{KP}. The only assumption about $\omega$ we need is
that each symbol $i$ has a limit frequency of appearance. For
simplicity, we will only present the proof for an example: a class
of iterated function schemes considered by Lalley and Gatzouras.

We refer the reader interested in other non-conformal random
iterated constructions to \cite{F}, \cite{GL2} and references
therein.

\section{Lalley-Gatzouras schemes}

The Lalley-Gatzouras scheme ${\cal F}$ is a self-affine IFS given
by a family of maps

\[
f_{i,j}(x,y) = (a_{ij} x + c_{ij}, b_i y + d_i),
\]
where the alphabet $A$ of allowed symbols is $1\leq i \leq m_1,
1\leq j\leq m_2(i)$. We will assume that $b_i\geq a_{ij}$ (that
is, the contraction in the horizontal direction is not weaker than
the contraction in the vertical direction for all maps).

We will also assume that for all $(i,j)\in A$ $0<a_{ij}<1, 0\leq
c_{i1}<\ldots <c_{im_2(i)}\leq 1-a_{im_2(i)}, c_{ij+1}\geq
a_{ij}+c_{ij}$ and that $0<b_i<1, 0\leq d_1<\ldots <d_{m_1}\leq
1-b_{m_1}, d_{i+1}\geq b_i+d_i$. We will say that the {\it
separation condition} holds if we actually have $c_{ij+1}>
a_{ij}+c_{ij}$ and $d_{i+1}> b_i+d_i$.

The main result of \cite{LG} is the formula for the Hausdorff
dimension of the limit set $\Lambda$:

\[
\dim_H(\Lambda) = \max \left\{ \frac {\sum_i \sum_j p_{ij} \log
p_{ij}} {\sum_i \sum_j p_{ij} \log a_{ij}} + \sum_i q_i \log q_i
\left( \frac 1 {\sum_i q_i\log b_i} - \frac 1 {\sum_i \sum_j
p_{ij} \log a_{ij}}\right)\right\},
\]
where $\{p_{ij}\}$ is a probability distribution on $A$,
$q_i=\sum_j p_{ij}$ and the maximum is over all possible
$\{p_{ij}\}$.

Consider now a family of Lalley-Gatzouras schemes $\{{\cal
F}_k\}_{k=1}^m$ with alphabets $A_k$ and maps $f_{i,j}^{(k)}$. As
mentioned above, we can apply them in any order $F_{\omega_1}\circ
F_{\omega_2}\circ\ldots$, $\omega=\omega_1 \omega_2\ldots \in
\Omega=\{1,\ldots,m\}^\N$ and obtain some limit set
$\Lambda_\omega$. We will assume that the limits

\begin{equation} \label{eqn:pi}
P_k = \lim_{n\to\infty} \frac 1 n \sharp \{1\leq l\leq n;
\omega_l=k\}
\end{equation}
exist and are positive. We will ask what is the value of
$\dim_H(\Lambda_\omega)$.

Before formulating the answer, let us note that any finite product
$F_{\omega_1}\circ\ldots\circ F_{\omega_n}$ is again a
Lalley-Gatzouras scheme. It follows that we can calculate the
Hausdorff dimension of $\Lambda_\omega$ for any periodic sequence
$\omega$. Given a {\it rational} probabilistic vector
$Q=(Q_1,\ldots,Q_m)$, we can choose a periodic sequence
$\omega(Q)$ in which the frequency of symbol $k$ is $Q_k$. Let us
write

\[
\L(Q)=\dim_H \Lambda_{\omega(Q)}.
\]
Our main result is as follows.

\begin{thm} \label{thm:main}
The function $\L(Q)$ is well defined, does not depend on the
choice of $\omega(Q)$. We can extend it by continuity to the whole
simplex of probabilistic vectors (we will keep the notation
$\L(Q)$ for the extended function). We have

\[
\dim_H(\Lambda_\omega) = \L(P).
\]
\end{thm}

\section{Proof of Theorem \ref{thm:main}}

Let us start by presenting a more detailed description of
$\Lambda_\omega$ (compare \cite{Hu}). Let
$A_\omega=A_{\omega_1}\times A_{\omega_2}\times\ldots$. We define
a projection $\pi_\omega:A_\omega\to\Lambda_\omega$ by the formula

\[
\pi_\omega((i_1, j_1), (i_2, j_2),\ldots) = \lim_{n\to\infty}
f_{i_1,j_1}^{(\omega_1)}\circ\ldots\circ
f_{i_n,j_n}^{(\omega_n)}(0,0),
\]
$(i_k,j_k)\in A_k$. We get

\[
\Lambda_\omega = \pi_\omega(A_\omega).
\]

Because of the nonconformality of the system, the most natural
class of subsets of $A_\omega$ to study are not cylinders but {\it
rectangles} (in particular, {\it approximate squares}). The
rectangle is defined as follows: given a sequence $(i,j)\in
A_\omega$ and two natural numbers $n_1\leq n_2$ we define

\[
R_{n_1, n_2}(i,j)=\{(i',j')\in A_\omega; i_k'=i_k \forall k\leq
n_2, j_k'=j_k \forall k\leq n_1\}.
\]

We will call

\[
d_1(R_{n_1, n_2}(i,j))=\prod_{k=1}^{n_1} a_{i_k j_k}^{(\omega_k)}
\]
the {\it width} and

\[
d_2(R_{n_1, n_2}(i,j))=\prod_{k=1}^{n_2} b_{i_k}^{(\omega_k)}
\]
the {\it height} of the rectangle $R_{n_1, n_2}(i,j)$. Indeed, the
projection of a rectangle under $\pi_\omega$ is the intersection
of $\Lambda_\omega$ with a geometric rectangle of the same width
and of the same height. The rectangle of approximately (up to a
constant) equal width and height is called an approximate square.

%We can introduce a natural metric on $A_\omega$. Given two
%sequences $(i,j), (i',j')\in A_\omega$ let $n_1((i,j),(i',j'))$ be
%the first $n$ for which $(i_n,j_n)\neq (i_n',j_n')$ and let
%$n_2((i,j),(i',j'))$ be the first $n$ for which $i_n\neq i_n'$. We
%write

%\[
%\rho_\omega((i,j),(i',j'))=\max
%(\rho_\omega^1((i,j),(i',j')),\rho_\omega^2((i,j),(i',j')))
%\]
%where

%\[
%\rho_\omega^1((i,j),(i',j')) = \prod_{k=1}^{n_1((i,j),(i',j'))-1}
%a_{i_k j_k}^{(\omega_k)},
%\]
%\[
%\rho_\omega^2((i,j),(i',j')) = \prod_{k=1}^{n_2((i,j),(i',j'))-1}
%b_{i_k}^{(\omega_k)}.
%\]
%It is easy to check that $\pi_\omega$ is a Lipschitz projection
%from $(A_\omega,\rho_\omega)$ to $\R^2$; if all the schemes ${\cal
%F}$ satisfy separation condition, $\pi_\omega$ is actually a
%bi-Lipschitz map and $\rho_\omega^1((i,j),(i',j'))$ and
%$\rho_\omega^2((i,j),(i',j'))$ are approximately the horizontal
%and vertical distances between $\pi_\omega(i,j)$ and
%$\pi_\omega(i',j')$.

Our main step is the following proposition.

\begin{prop} \label{en}
For $Q$ a rational probabilistic vector close to $P$ and for any
choice of $\omega(Q)$, we can construct a bijection
$\tau:A_{\omega(Q)}\to A_{\omega}$ with the following properties.
Let $R=R_{n_1, n_2}^{(\omega(Q))}(i,j)$ be an approximate square
in $A_{\omega(Q)}$ of width $d$. If $\delta = \max |P_k-Q_k|$ is
sufficiently small, $\tau(R)$ contains an approximate square of
width at least $d^{1+K\delta+\epsilon}$ and is contained in an
approximate square of width at most $d^{1-K\delta-\epsilon}$, $K$
depending only on the iterated schemes but not on $P$ or $Q$ and
$\epsilon$ arbitrarily small for sufficiently small $d$.
\begin{proof}

We will need the following simple statement (a reformulation of
\eqref{eqn:pi}):

\begin{lem} \label{random}
For every $n$ there exists $\epsilon(n)$ such that for each $k$
the $n$-th appearance of symbol $k$ in the sequence $\omega$ takes
place between positions $n/P_k(1-\epsilon(n/P_k))$ and
$n/P_k(1+\epsilon(n/P_k))$. Moreover, $\epsilon(n)$ goes
monotonically to 0 as $n$ goes to $\infty$.
\end{lem}

Consider now the pair of sequences: $\omega$, the sequence we work
with, and $\omega(Q)$, a periodic sequence with frequencies $Q$.
We will assume that $Q$ is $\delta$-close to $P$ and that both
probabilistic vectors are positive. Obviously, in the sequence
$\omega(Q)$ the $n$-th appearance of symbol $k$ is at position
$n/Q_k$, give or take a constant.

We will define $\chi_{\omega, \omega(Q)}$ as a permutation of $\N$
in the following way: if $l_1$ is the place of $n$-th appearance
of symbol $k$ in the sequence $\omega$ and $l_2$ is the place of
$n$-th appearance of symbol $k$ in the sequence $\omega(Q)$, we
set $\chi_{\omega, \omega(Q)}(l_1)=l_2$. We can then construct a
bijection $\tau:A_{\omega(Q)}\to A_{\omega}$ as

\[
\tau((i_1,j_1), (i_2, j_2),\ldots) = (i_{\chi_{\omega,
\omega(Q)}(1)},j_{\chi_{\omega, \omega(Q)}(1)}) \ldots
\]

Denote

\[
D_1 = \chi_{\omega, \omega(Q)}(\{1,\ldots,n_1\})
\]
and
\[
D_2 = \chi_{\omega, \omega(Q)}(\{n_1+1,\ldots,n_2\})
\]

We remind that the rectangle $R$ is defined as the set of
sequences $(i',j')\in A_{\omega(Q)}$ for which we fix the first
$n_1$ $(i_k',j_k')$ and the following $n_2-n_1$ $i_k'$. Hence, the
set $\tau(R)$ is the set of sequences $(i',j')\in A_{\omega}$ for
which we fix $(i_k',j_k')$ for $k\in D_1$ and we fix $i_k'$ for
$k\in D_2$.

Denote
\[
r_1 = \inf(\N\setminus D_1) -1,
\]
\[
r_2 = \inf(\N\setminus (D_1\cup D_2)) -1,
\]
\[
s_1 = \sup(D_1),
\]
\[
s_2 = \sup(D_1\cup D_2).
\]
We have

\[
R_{s_1,s_2}^{(\omega)}(\tau(i,j)) \subset \tau(R) \subset
R_{r_1,r_2}^{(\omega)}(\tau(i,j)).
\]

Assume $\delta$ is much smaller than any $P_k$. By Lemma
\ref{random},

\[
r_1 \geq n_1(1-K_0\epsilon(n_1)-K_0\delta),
\]
\[
r_2 \geq n_2(1-K_0\epsilon(n_2)-K_0\delta),
\]
\[
s_1 \leq n_1(1+K_0\epsilon(n_1)+K_0\delta),
\]
\[
s_2 \leq n_2(1+K_0\epsilon(n_2)+K_0\delta)
\]
for some $K_0>0$ depending only on the iterated schemes.

Consider the width of $R_{r_1,r_2}^{(\omega)}(\tau(i,j))$ versus
the width of $R$. The latter is a product of $n_1$ numbers
$a_{ij}^{(k)}$, the former it the subproduct of $r_1$ of those
numbers. As all $a_{ij}^{(k)}$ are uniformly bounded away from 0
and 1,

\[
d_1(R_{r_1,r_2}^{(\omega)}(\tau(i,j))) \leq
d^{1-K\epsilon(n_1)-K\delta}
\]
for some uniformly chosen $K$, depending only on the iterated
schemes. Similar reasoning proves

\[
d_2(R_{r_1,r_2}^{(\omega)}(\tau(i,j))) \leq
d^{1-K\epsilon(n_2)-K\delta}.
\]

Consider now the width of $R_{s_1,s_2}^{(\omega)}(\tau(i,j))$
versus the width of $R$. The former a product of $s_1$ numbers
$a_{ij}^{(k)}$, the latter it the subproduct of $n_1$ of those
numbers, the same reasoning as before gives us

\[
d_1(R_{s_1,s_2}^{(\omega)}(\tau(i,j))) \geq
d^{1+K\epsilon(n_1)+K\delta},
\]

\[
d_2(R_{s_1,s_2}^{(\omega)}(\tau(i,j))) \geq
d^{1+K\epsilon(n_2)+K\delta}.
\]
The rectangles $R_{r_1,r_2}^{(\omega)}(\tau(i,j))$ and
$R_{s_1,s_2}^{(\omega)}(\tau(i,j))$ are not necessarily
approximate squares, but we can easily replace the former by some
slightly larger rectangle which is an approximate square and we
can replace the latter by some slightly smaller rectangle which is
an approximate square. We are done.
\end{proof}
\end{prop}

\begin{rem}
We can introduce a metric on $A_\omega$, defining the distance
between two points as the sum of width and height of the smallest
rectangle containing them both. This metric is natural because if
the maps satisfy separation condition, $\pi_\omega$ is
bi-Lipschitz (without separation condition it will only be a
Lipschitz projection). In this metric, the maps $\tau$,
$\tau^{-1}$ are H\"older continuous with every exponent smaller
than 1 (if $\delta=0$) or with exponent $1-K\delta$ (if $\delta$
is positive but small).
\end{rem}

This proposition basically ends the proof of Theorem
\ref{thm:main}. By Proposition 3.3 and Lemma 5.2 in \cite{LG}, for
any Lalley-Gatzouras scheme there exists a probabilistic measure
$\mu$ supported on $A^\N$ such that
\begin{itemize}
\item[i)] for a $\mu$-typical point $(i,j)$ and the decreasing
sequence of all approximate squares $R_k=R_{n_1(k), n_2(k)}(i,j)$,
\[
\frac {\log \mu(R_k)} {\log d_1(R_k)} \to \dim(\Lambda),
\]

\item[ii)] for every point $x\in A^\N$ there exists a decreasing
sequence of approximate squares $R_k=R_{n_1(k), n_2(k)}(i,j)$ for
which

\[
\frac {\log \mu(R_k)} {\log d_1(R_k)} \to \dim(\Lambda).
\]
\end{itemize}

We can define such measure $\mu_Q$ supported on $A_{\omega(Q)}$
for any rational $Q$ (because this is again a Lalley-Gatzouras
scheme). We can then transport this measure to $A_\omega$ by the
map $\tau$. We obtain a measure $\nu_Q$ such that

\begin{itemize}
\item[i)] for a $\nu_Q$-typical point $(i,j)$ and the decreasing
sequence of all approximate squares $R_k=R_{n_1(k), n_2(k)}(i,j)$,
\[
\liminf \frac {\log \nu_Q(R_k)} {\log d_1(R_k)} \geq \L(Q)
(1-K\delta),
\]

\item[ii)] for every point $x\in A_\omega$ there exists a
decreasing sequence of approximate squares $R_k=R_{n_1(k),
n_2(k)}(i,j)$ for which

\[
\limsup \frac {\log \nu_Q(R_k)} {\log d_1(R_k)} \leq \L(Q)
(1+K\delta).
\]
\end{itemize}

It implies that
\[
\L(Q) (1-K\delta) \leq \dim_H \Lambda_\omega \leq \L(Q)
(1+K\delta),
\]
the proof is as in \cite{LG}. \qed

This result has immediate applications for random systems,
obtained by choosing $\omega$ randomly with respect to some
Bernoulli measure on $\Omega$.

On the other hand, this method is not going to work for
stochastically-selfsimilar systems considered in \cite{F} or
\cite{GaL}. For such systems we would not have a single sequence
$\omega$ but instead $\omega$ would depend on the point in the
fractal. While we would still be able to define $\tau$ almost
everywhere, the sequences $\omega(x)$ at different points
$x\in\Lambda$ would not all satisfy Lemma \ref{random}, and hence
$\tau$ would not everywhere have nice H\"older properties.

%\newpage
\bibliography{ref}

\end{document}